# Advanced Quantitative Techniques to Solve Center of Gravity Problem in Supply Chain

**Brian Houck**
*PwC US*
brian.m.houck@pwc.com

**Chetan Sampat**
*PwC US*
chetan.p.sampat@pwc.com

**Srijit Maiti**
*PwC US*
srijit.maiti@pwc.com

**Shivam S**
*PwC US*
shivam.s.s@pwc.com

**Anurag Vaishistha**
*PwC US*
anurag.m.vashistha@pwc.com

**Sumit Banerjee**
*PwC US*
sumit.b.banerjee@pwc.com

*Abstract*—Activities involving transformation of raw materials, various resources and components into final products and also delivering it to the end customer incur a significant cost during the selection of location of a warehouse that can be easily accessed by various actors of the supply chain. To minimize upstream and downstream transportation costs, the center of gravity (CoG) analysis method is used to find the potential warehouse locations for a given demand network which have an impact on the entire supply chain network. Mixed Integer Linear Programming (MILP), an open source tool is developed for implementing CoG method along with certain service level constraints to find optimal potential locations with the least cost. In this paper, an optimization tool has been designed for a forward logistics network with several novel methods like "Customer Location Selection" (CLS), "Customer Packets" along with other business heuristics that optimize and enhance the existing MILP to get the optimal solutions with low computational cost and runtime. Finally, recommending an alternative network of facilities which reduces overall costs compared to the existing network. An user interface has also been developed to make a user friendly interaction with the model. We can conclude that this model can significantly help companies reduce costs during the logistics network design.

## I. INTRODUCTION

Warehousing is a key driver of the supply chain performance, providing economies of scale through efficient operations, storage capacity and a central location with the goal of minimizing cost and fulfillment times. In today's world, the selection of a warehouse location has become one of the most imperative and strategic decisions in the optimization of logistic systems. Warehouse location decision is one of the strategic decisions in the fields of supply chain management (SCM) and operations management (OM). The moot point is that as it requires large investment and since it is irrecoverable in most of the instances, choice of a warehouse location among eligible alternatives is a very powerful decision. A poor choice of location might result in excessive transportation cost, loss of qualified labor, competitive advantage or some similar condition that would be detrimental to operations. As decisions regarding facilities are a crucial element in a company's quality success or failure, the company must analyze in order to find a balance between responsiveness and efficiency that best matches the needs of its target customer. Warehouses are key aspects of modern supply chains and play a vital role in the success, or failure of businesses today.

In order to reduce transportation cost, enforce operation efficiency and logistic performance, evaluating and selecting a suitable Distribution Center (DC) location has become one of the most important decision issues for distribution industries. In the process of selection it is necessary first and foremost to identify the set of influential factors relevant to the DC location selection. Many influential factors are considered for the selection of a particular plant location, e.g. investment cost, climate condition, labor force quality and quantity, transportation availability, etc. [1,2,3,4,5]

The problem of warehouse location selection contains an efficient strategic investment decision on long term and business profitability. At this point, one of the most important decision making processes of logistics administrators is the location decision of the distribution center. [6]. In addition to the cost and profitability indicators under increasing competition conditions of today customer service level and customer satisfaction are also demanded. [7]

The warehouse location selection is a process of selecting an allocation center in an economic region where there are some supply stations and certain demand points. A typical facility location problem consists of choosing the best among potential sites, subject to constraints requiring that demands at several points must be serviced by the established facilities. The objective of the problem is to select facility sites in order to minimize costs; these typically include a part which is proportional to the sum of the distances from the demand points to the servicing facilities, in addition to costs of opening them at the chosen sites.

The classical algorithms include the Center of Gravity (CoG) approach, capacitated facility location problem model, and the P-median selecting location model, in the warehouse location selection. The Center of Gravity (CoG) Method is an approach that seeks to compute geographic coordinates for a potential single new facility that will minimize costs. It's an approach where the main inputs that it considers are Markets, Volume of goods shipped, Shipping costs. This method is beneficial because of it's simplified computation, consideration of existing facilities, and it minimizes costs. This model uses several assumptions, namely: 1) The cost of transportation costs is assumed to increase in proportion to the volume transferred. 2) Good supply sources and the location of the production can be located on a map with clear X and Y coordinates. The accuracy of computation is low, hence the result can only be referenced. [8]

The capacitated facility location problem is the basis for many practical optimization problems, where the total demand that each facility may satisfy is limited, i.e. facilities have constraining capacities on the amount of demand they can serve. Hence, modeling such problems must take into account both demand satisfaction and capacity constraints. [9] The complex computation process is the shortcoming of the model.

A median problem states that given a number of facilities from possible points in a graph, in such a way that the sum of the distances from each customer to the closest facility is minimized. Often, the number of facilities to be selected is predetermined in advance; this number is commonly denoted p, and the corresponding problem is given this symbol as a suffix. The p-median problem is therefore a variant of the uncapacitated facility location problem and specifically seeks to establish p facilities, without considering fixed costs.

In this paper, we show PwC's approach using various catalyzing mechanisms of business intelligence to speed up the traditional CoG methods which resulted in remarkable improvements in the capabilities of MIP algorithms. Some of the biggest contributors have been pre-solve, cutting planes and using heuristics. Presolving the MIP includes our novel CLS method and Clustering techniques which helps to confine our warehouses to more specific regions. Cutting planes uses techniques to tighten the formulation by removing undesirable solutions using step well optimization of first finding suboptimal solutions to effectively scope out the problem, basically breaking down the initial big problem into two steps: Country level and State level. Along with these, we have also used certain heuristics to help the MIP solver by constraining the flow variables to be of Boolean type and creating customer packets. Using the aforementioned techniques, we have developed a tool with an easy to use Excel interface that is designed to be used with little or no prior training.

## II. APPROACH

A. Problem Definition

CoG Analysis is a quick way to identify optimal distribution center (DCs) locations for a given demand network answering the following important questions:
1. Where should distribution centers be geographically located to minimize cost?
2. Which customers will be supplied from each distribution center?
3. To visualize the distribution and evaluate Scenarios
4. Conduct demand sensitivity tests
5. Simulate Future network and check resilience of the network

In this paper, we aim to solve these key business problems using Mixed Integer Linear Programming (MILP) and develop some novel techniques to catalyze the optimisation process.

B. Qualitative formulation of problem statement
- Finding a probable location for a Warehouse: We aim to solve these two questions here:
    a. To determine the best geographic locations of potential "greenfield" sites given a set of customers.

b. To know how would the footprint look like by determining the Greenfield and how much volume each Greenfield moves
- Building an Incremental Network: We are predicting the following:
   a. To predict the new Greenfield location and the footprint after keeping "X" number of facilities in the current network

- Service level Analysis: Following advanced analysis is being carried out with number of additional Service level constraints
   a. Average weighted Distance of Greenfield should be less than "Z"
   b. To serve "X" % of demand with "Y" miles from a Greenfield

All of these three problems constructed qualitatively above would work as constraining the main transportation objective cost function which is formulated in the next section.

C. Mathematical Formulation of Objective Function

Objective function is the function of cost which needs to be minimized and in this case, its transportation cost. Before defining the objective function, let us understand the decision variables that we will be solving for. First decision variable is c which denotes whether a warehouse should be opened at that location or not. It is a binary variable which when 0 denotes no warehouse is needed at that location and 1 denotes opening of new warehouse or existing warehouse at that location. The other decision variable that we need to find out by optimizing the cost function is flow between the opened or existing warehouse and consumer. If there are y probable warehouses and x customer demand points, the total decision variables would be y for warehouse decision and x*y for flow decision.

Cost Function $= \sum f_{ij} * d_{ij} + \sum c_j$

Subjected to constraints:

- For decision variables initialisation

$f_{ij} \geq 0$

$c_j = \{0, 1\}$

- Warehouse constraint

$\sum_{j=1}^{y} c_j =$ Warehouse limit

- Demand constraints

$\sum_{j=1}^{y} f_{ij} \geq$ Demand of $C_i$ customer

- Conservation of flow

$\sum_{i=1}^{x} \sum_{j=1}^{y} f_{ij} \geq$ Total demand of all customers

- Maximum average distance constraint

$\sum (f_{ij} * d_{ij})/Total\ demand \leq MAD$

- Percentage of households served within constraint

$\sum (f_{ij} * p_{ij})/Total\ demand \geq MPCT$

Notations used:
$f_{ij}$ : Flow or weight at i-j combination of customer-warehouse
$d_{ij}$ : Distance from location customer i to warehouse j
$ij$ : the ij[th] combination of customer-warehouse arc
$j$ : used for notation of jth warehouse
$i$ : used for notation of ith customer
x : Number of customer demand points
y : Number of probable warehouse locations
*MAD* : Maximum average distance
*MPCT* : Minimum percentage of households served within a particular distance from warehouse

D. Methodology
*'Data gives you the what, but humans know the why'*

Four novel methods combined with traditional MILP using business intelligence have been developed to enhance the existing branch and bound technique and decrease the running time for optimizing objective function which is a NP hard problem.
MILP problems are generally solved using a linear-programming based *branch and bound algorithms*. Basic LP-based branch-and-bound can be described as follows. We begin with the original MIP. Not knowing how to solve this problem directly, we remove all of the integrality restrictions. The resulting LP is called the linear-programming relaxation of the original MIP. We can then solve this LP. If the result happens to satisfy all of the integrality restrictions, even though these were not explicitly imposed, then we have been quite lucky. This solution is an optimal solution of the original MIP, and we can stop. If not, as is usually the case, then the normal procedure is to pick some variable that is restricted to be integer, but whose value in the LP relaxation is fractional.

For the sake of argument, suppose that this variable is x and its value in the LP relaxation is 5.7. We can then exclude this value by, in turn, imposing the restrictions x ≤ 5.0 and x ≥ 6.0. If the original MIP is denoted P0, then we might denote these two new MIPs by P1, where x ≤ 5.0 is imposed, and P2, where x ≥ 6.0 is imposed. The variable x is then called a branching variable, and we are said to have branched on x, producing the two sub-MIPs P1 and P2. It should be clear that if we can compute optimal solutions for each of P1 and P2, then we can take the better of these two solutions and it will be optimal to the original problem, P0. In this way we have replaced P0 by two simpler (or at least more-restricted) MIPs. We now apply the same idea to these two MIPs, solving the corresponding LP relaxations and, if necessary, selecting branching variables. In doing so we generate what is called a search tree.

The MIPs generated by the search procedure are called the nodes of the tree, with P0 designated as the root node. The leaves of the tree are all the nodes from which we have not yet branched. In general, if we reach a point at which we can solve or otherwise dispose of all leaf nodes, then we will have solved the original MIP.

Now, let us discuss the novel methods that can be used to enhance the traditional MILP

a. Step-Well Optimization

An Objective function, if left to optimize on its own using a set of constraints, would consider each possible combination, calculate the cost and would provide us with the solution with optimized objective value in a horizontal fashion. This would definitely provide us the optimality but at the same time, comes with a huge computational cost and runtime. For COG Analysis, we want optimality but also low computational cost and runtime. To solve this problem, we introduced a top down approach to optimization. Instead of finding the optimal local town and street where the warehouse is supposed to be opened at the very start, we would be finding the suboptimal points which would represent the optimal state where the warehouse needs to be opened and then the optimal point or the proper location of the warehouse in that state could be found. This led to major reduction in the run time by just adding a layer to the problem.

Mathematically, we take the gradient of the cost function and check how the cost function changes. As, if the gradient or the slope of the Cost function is proportional to higher orders of the latitude and longitude then we can't simply apply this step well method and we have to check the cost function at all points instead.

$$\nabla Cost\ function = \nabla(\sum_{i \in I} f_i * \sqrt{(x - x_i)^2 + (y - y_i)^2})$$

$$= \partial/\partial x(\sum_{i \in I} f_i * \sqrt{(x - x_i)^2 + (y - y_i)^2})\ \hat{i}$$

$$+ \partial/\partial y(\sum_{i \in I} f_i * \sqrt{(x - x_i)^2 + (y - y_i)^2})\ \hat{j}$$

$$= (\sum_{i \in I} f_i * (x - x_i)/(\sqrt{(x - x_i)^2 + (y - y_i)^2}))\ \hat{i} +$$

$$(\sum_{i \in I} f_i * (y - y_i)/(\sqrt{(x - x_i)^2 + (y - y_i)^2}))\ \hat{j}$$

As we can see that the order of the gradient is low and consequently the curve is smooth and not volatile, we can use this step well approach to calculate cost at selected points.

In order to select a set of initial warehouse locations from probable warehouse locations that are available to us as an input, we used Candidate Location Selection (CLS) Analysis.

b. Candidate Location Selection (CLS) Analysis

We have developed a novel method called Candidate Location Selection (CLS) analysis to reduce the size of a problem by selecting the best possible warehouse candidate locations for the suboptimal solution. This technique can be typically applied in advance of the start of the branch-and-bound procedure. These reductions are intended to reduce the size of the problem and to tighten its formulation. Using the inherent relationships present in the data, we can select how many candidate locations are needed to represent the entire state based on three main pillars.

- **Area**: Bigger the State, more should be the representation via warehouse locations. Represented by Ascore
- **Proximity to the existing warehouse**: Closer the existing warehouses, less would be the chances of opening a warehouse nearby in that state. Represented using Pscore
- **Demand Density of the State**: Higher the demand density in the state, higher the chances of having a warehouse there. Represented by Dscore

All these scores calculated are not absolute but relative and are scaled between 1-10. i.e If the area of state S is smallest, than it's Ascore will be 1. Using these three metrics, we calculate an overall average score for the state using which we assign the number of input probable warehouse locations for that state.

$$CLS\ Score\ =\ \sqrt[3]{Ascore\ *\ Pscore\ *\ Dscore}$$

Now, depending upon CLS Score, we provide the number of locations to each state that would represent the state in the initial input for optimization. So, as suggested in the step well approach above, we don't need all the warehouse inputs in the first step itself. Hence, CLS technique is used in the first step to get the probable initial inputs that would be used to calculate the suboptimal locations of warehouses

c. Clustering

Clustering can be used as a presolve and to get an initial probable warehouse candidate location for suboptimal solutions when there is no existing warehouse present. Clustering intuitively reduces the transportation cost when the distances are weighted using demand flow on it.
It is a Centroid-based clustering, the central vector represents the number of clusters which is not necessary to be a member of the data set. When we fix the number of clusters to k, a formal definition of an optimization problem is given by k-means clustering: finding the k cluster centers and assign the objects to their nearest cluster center in the following manner to the cluster which gets the minimum squared distance.

d. Customer Packets

Customer Packets are the group of closely located Customers that would have high probability of receiving the orders from the same warehouse. So, we can reduce the dimension of the problem to a large extent by using this heuristic to our benefit in Branch and Bound.
That means instead of using for example 900 Customer demand points, we can use these newly created hypothetical 150 Customer Packets and their demand and distance from the warehouse candidate locations. Customer Packets are formed in a way that Overall Objective function is not affected by the new coordinates or new data points of the packets as Distance metric is created using the following formula.
Suppose for example, our Cost function looks something like:

$$f_{11}*d_{11} + (f_{12}*d_{12} + f_{13}*d_{13} + f_{14}*d_{14}) + f_{15}*d_{15} + f_{16}*d_{16} + f_{17}*d_{17} + f_{18}*d_{18} + \ldots$$

Now, if $C_2\ C_3\ C_4$ are located very close by with respect to the warehouse then we can treat them as a packet, and can resolve the equation to

$$f_{11}*d_{11} + f_{1p}*d_{1p} + f_{15}*d_{15} + f_{16}*d_{16} + f_{17}*d_{17} + f_{18}*d_{18} + \ldots$$

where $f_{1p} = f_{12} + f_{13} + f_{14}$
and $d_{1p} = (f_{12}*d_{12} + f_{13}*d_{13} + f_{14}*d_{14})/f_{1p}$
Hence, reducing the number of decision variables which were earlier three ($f_{12}, f_{13}, f_{14}$) to only one ($f_{1p}$) which helps reduce the computational cost for finding the optimal location.

## III. RESULTS AND INSIGHTS

The major concerns for an optimization problem or a NP hard problem are the huge costs and run time associated with them, for which we developed several novel approaches as discussed in the approach section to minimize the run time along with maintaining high accuracy and stable results.
Run time of Center of Gravity optimization problem depends upon three major external factors as shown in the figure below

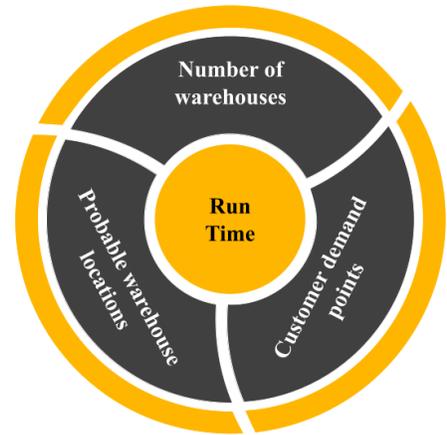

Top to Bottom Optimization along with CLS and Customer packets focuses on solving Warehouse and customer demand points problem respectively and hence, making run time of the model independent of them.
In a simple MILP, if Customer Count or Warehouse Count is n which when increased x times the run time increases by $8^{(x-1)}$ approximately. There are other internal factors also that could affect the run time like the type of unknown variables in the Cost function, their limits, separation of warehouse domain selection along with other heuristics that we were able to optimize using our cutting domain methodology and developing some new heuristics.

We have used demand data along with coordinates of a preferred time period of the user and existing warehouse location coordinates. The data fields required for the computation were: Demand, Demand Latitude, Demand Longitude, Existing Warehouse State Name, Existing Warehouse Latitude, Existing Warehouse Longitude, Existing Warehouse Operation Status.

A. CLS Analysis results: Used for Warehouse Count problem

For each state, we utilized CLS analysis to solve for an optimal number of probable warehouse locations which allowed us to see what is the distribution of probable warehouses in the states. Next, we applied clustering to locate the warehouse locations in each state. We used the COG facility scenarios to decide how many facility nodes per state we would have the MILP select. This led to an overall decrease in run time and independence of the Run time wrt warehouse count as can be seen in the above figure. Orange bar which represents the simple MILP without top bottom approach and CLS analysis is highly dependent on number of probable warehouse inputs, while blue bar representing pwc approach can be seen as independent of increasing warehouse count.

Table 1: Representative CLS score & whether it will be located on ridge or trough

| State | AScore | PScore | DScore | Allocation |
|---|---|---|---|---|
| California | 1 | 1 | 1 | Trough (1) |
| New york | 1 | 1 | 0 | Near to trough(2) |
| Penn | 1 | 0 | 0 | Near to Ridge |
| Vermont | 0 | 1 | 1 | Near to Trough |
| Florida | 0 | 0 | 1 | Near to Ridge(3) |
| Washington DC | 0 | 0 | 0 | Ridge(4) |

Figure 1: Cost associated with the States in table above

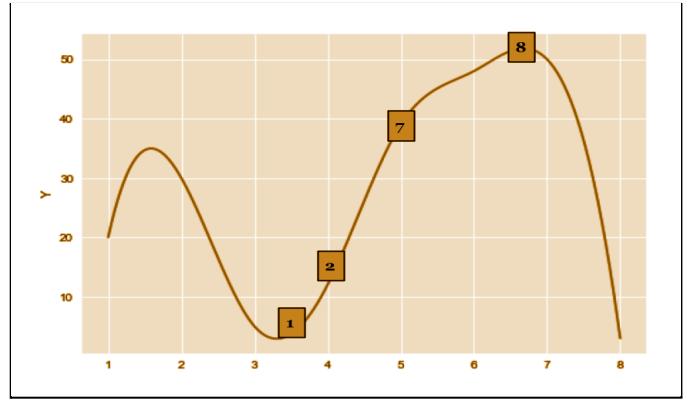

B. Customer Packets Analysis: Used for Customer Count problem

Customer Packets were made when there were multiple points located very close to one another and the entire packet was treated as a single demand point which led to an overall decrease in count of Customers. As can be seen in the figure above, Increase in customer count leads to increase in Run time for simple MILP but not for Modified MILP with Customer packets formation.

C. MILP Model Results: Combining all of the approaches to enhance the optimization

The MILP model generated the cost information from the inputs we provided. We ran multiple scenarios for each of the states, the results summarized here are for the 1 and 2 COGs, where we forced the MILP to use greenfield locations and the company's existing facilities.

The first observation we made from the MILP results was related to the transportation costs. We wanted to compare the difference in weighted average distance (miles) of warehouse locations in the new network compared to an existing tool. The overall difference was 18 miles.

Table 2: PwC results and comparison with Benchmark

| Comparison Analysis | Benchmark Results | | PwC Results | | |
|---|---|---|---|---|---|
| Analysis Type | Latitude | Longitude | Latitude | Longitude | Diff miles |
| 1 COG | 40.4011 | -78.4519 | 40.4831 | -78.3519 | 17.2 |
| 1 COG | 37.85692 | -87.5325 | 37.540 | -87.3148 | 9.2 |

| 1 COG | 40.59428 | -75.6487 | 40.6018 | -75.5032 | 8.6 |
| 2 COG | 35.71492 | -118.815 | 35.8349 | -119.0902 | 25.4 |
| 2 COG | 40.59428 | -75.6487 | 40.5092 | -76.0434 | 20.1 |
| 2 COG | 35.93612 | -117.905 | 36.0121 | -117.7892 | 20.1 |
| 2 COG | 38.18434 | -83.2804 | 38.2721 | -83.0899 | 17.2 |

Figure 2: Runtime comparison with & without novel techniques

a. With CLS Analysis

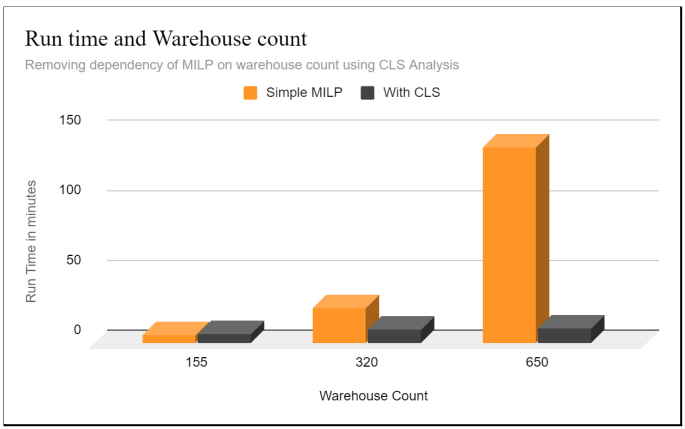

b. With Customer Packets

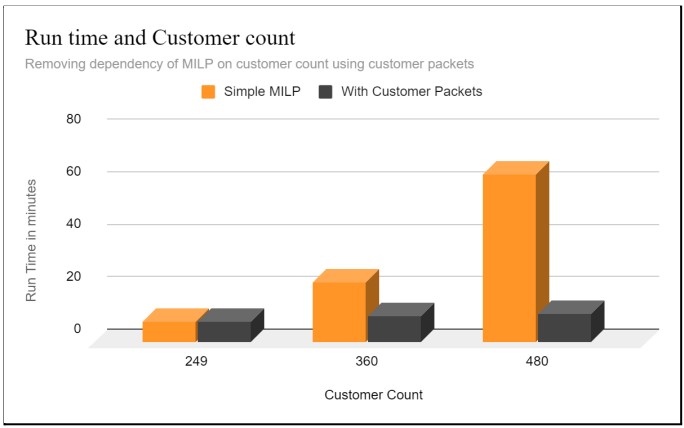

## IV. CONCLUSION

This paper proposes 4 different novel quantitative techniques to perform center of gravity analysis in supply chain design. The techniques have been proven to demonstrate accurate results as compared to those of existing tools which perform similar analysis (18 miles of difference). The methodology proposed is useful for any greenfield analysis. The first three techniques namely Step Well Optimization, CLS Analysis and Clustering would help us solve the problem when the size of input count for probable warehouses is huge, while Customer Packets could be used to increase the efficiency and to bring down the run time when size of the input demand points is huge. As future work, an improvement of model formulation will be explored by adding various additional parameters and constraints like capacitated warehouses, optimal number of warehouses needed etc.